\input amssym.def
\input  amssym.tex

\documentstyle[11pt]{article}

\title{The Consistency of $ZFC+CIFS$}

\author{Garvin Melles\thanks{Would like to thank Ehud Hrushovski
for supporting him with funds from NSF Grant DMS 8959511}\\Abraham
Fraenkel Center for Mathematical Logic 
\\Bar Ilan University\\Institute of Mathematics\\Hebrew University of Jerusalem}

\newcommand{\proof}{{\sc proof} \hspace{0.1in}}

\newcommand{\sub}{\subseteq}
\newcommand{\al}{\alpha}
\newcommand{\be}{\beta}
\newcommand{\ga}{\gamma}

\begin{document}
\mathsurround=.1cm
\maketitle

\begin{center}
{\bf INTRODUCTION}
\end{center}

This paper is a technical continuation of ``Natural Axiom Schemata
Extending ZFC. Truth in the Universe?'' In that paper we argue that
$CIFS$ is a natural axiom schema for the universe of
sets. In particular it is a natural closure condition on $V$ and a natural
generalization of $IFS(L).$ Here we shall prove the  
consistency of $ZFC\ +\ CIFS$ relative to the
existence of a transitive model of $ZFC$ using the compactness theorem together
with a class forcing.

\begin{center}
{\bf GENERAL FRAMEWORK FOR THE CLASS FORCING}
\end{center}

\noindent {\bf Notation-} $\exists !P\psi(P)$ will mean there is a
unique separative partially ordered set $P$ with maximal element such
that $\psi(P).$  

\vspace{.1in}

\noindent {\bf Definition 1.} \ Let $\psi(x,y)$ be a formula (in general
with a hidden parameter.) A partial order ${\bf P}$ (a proper class)
with order ${\bf \leq}$ is
said to be a $\psi(x,y)$ definable $Ord$ iteration if
$${\bf P}=\bigcup\limits_{\alpha\,\in\, Ord-\{0\}}P_{\alpha}$$
and
$${\bf \leq}\  =\bigcup\limits_{\alpha\,\in\,Ord-\{0\}}\leq_{\alpha}$$ 
where each
$P_{\alpha}$ is a set of $\alpha$ sequences and $P_{\alpha}$ and
$\leq_{\alpha}$ are defined by induction on $\alpha$ as follows:\\ 
$$P_1=\big\{(0,p)\mid p=1\ \vee\ V\models \exists !P(\psi(P,0)\ \wedge\
p\in P)\big\}$$
If $\alpha$ is a limit then $p\in P_{\al}$ iff $p$ is an $\alpha$
sequence and $\forall \beta<\alpha
\ (\beta\neq 0)\ p\restriction \beta\in P_{\beta}.$ For $p,q\in
P_{\al},$ $p\leq_{\alpha} q$ iff
$\forall\beta<\alpha\ (\beta\neq 0) \ p\restriction
\beta\leq_{\beta}q\restriction \beta.$\\

\vspace{.15in}

\noindent $P_{\alpha+1}=\Big\{p\frown \dot{p}\mid p\in P_{\alpha}\ \wedge\ \dot{p}=1\ \vee$

\begin{flushright}
$\dot{p}\hbox{ is
a }P_{\alpha}\hbox{ name } \wedge\ p\Vdash \exists !P(\psi(P,\alpha)\
\wedge\ \dot{p}\in P)\Big\}$
\end{flushright}

\vspace{.15in}

\noindent $\leq_{\alpha+1}\ =\ \Big\{(p\frown\dot{p},q\frown\dot{q})\mid p\leq_{\alpha}q\ \wedge\
\dot{q}=1\ \vee$

\begin{flushright}
$p\Vdash \exists !P(\psi(P,\alpha)\ \wedge\
\dot{p}\in P\ \wedge\ \dot{q}\in P\ \wedge\ \dot{p}\leq \dot{q}\,)\Big\}$
\end{flushright}

\vspace{.15in}

\noindent {\bf Remark 1.} \ The elements of ${\bf P}$ are really equivalence
classes induced by the relation $x\sim y\ \leftrightarrow\ x\leq y\
\wedge\ y\leq x.$

\vspace{.1in}

\noindent {\bf Remark 2.} \ If $p\in P_{\alpha}$ such that 
$$p\Vdash \exists !P\psi(P,\alpha)$$
then we identify $p\frown 1$ and $p\frown\dot{p}$ where $p\Vdash \dot{p}$ is
the maximal element of $P.$ 

\vspace{.1in}

\noindent {\bf Remark 3.} \ If $p\in {\bf P}$ then we identify $p$ and
$p\frown \bar{1}$ where $\bar{1}$ is any sequence of ones. 

\vspace{.1in}

\noindent{\bf Lemma 1.} \ Let ${\bf P}$ be a $\psi(x,y)$ definable
iteration. Then ${\bf P}$ and each $P_{\alpha}$ are separative. \\
\proof The case $\alpha$ is a limit ordinal we leave to the reader. 
By the definition of $\psi(x,y)$ definable iteration $P_1$ is
separative. Let $(p,\dot{p})$ and $(q,\dot{q})\in P_{\alpha+1}$ such
that $(p,\dot{p})\not\leq (q,\dot{q}).$ If $p\not\leq q,$ then we are
done by the induction hypothesis. So suppose $p\leq q.$ Then 
$$p\not\Vdash \dot{p}\leq \dot{q}$$
so there is an $r\leq p$ such that
$$r\Vdash \dot{p}\not\leq \dot{q}$$
so for some $\dot{r}$
$$r\Vdash \dot{r}\leq \dot{p}\ \wedge\ \dot{q}\perp \dot{r}$$
$(r,\dot{r})\leq_{\alpha+1}(p,\dot{p})\ \wedge\ (r,\dot{r})\perp (q,\dot{q}).$

\vspace{.1in}

\noindent{\bf Definition 2.} \ Let ${\bf P}$ be a $\psi(x,y)$ definable
iteration. Let
$\alpha\in Ord,$ and let $G$ be a generic subset of $P_{\al}.$ Let
$\varphi(x_1,\ldots,x_n)$ be a formula and $a_1\ldots, a_n\in
V^{P_{\alpha}}.$ For every $\alpha<\beta\in Ord,$ we define in $V[G]$ by
induction on $\be$ maps
$\pi_{\alpha\beta},\ \pi'_{\alpha\beta},\ \pi''_{\alpha\beta}$ and
partial orders $(P_{\al\be},\ \leq_{\al\be})$ as follows:

\noindent Let $\beta=\alpha+1.$ If $(p,\dot{q})\in P_{\alpha+1},$
and $p\in G,$ then 
$\pi_{\alpha\,\alpha+1}(p,\dot{q})=1$ if $\dot{q}=1$ and
$$\pi_{\alpha\,\alpha+1}(p,\dot{q})=i_G(\dot{q})$$
otherwise. If $p\not\in
G,$ then $\pi_{\al\,\al+1}(p,\dot{q})$ is not defined.  Let
$$P_{\al\,\al+1}=\big\{\pi_{\al\,\al+1}(p)\mid p\in P_{\al+1}\}$$ 
and let
$$\leq_{\al\,\al+1}=\big\{(\,\pi_{\al\,\al+1}(p),\ \pi_{\al\,\al+1}(q)\,)\mid
(p,q) \in\  \leq_{\al+1}\}$$

\noindent If $U\sub P_{\al+1}$ is a
regular cut, then
$$\pi'_{\al\,\al+1}(U)=\overline{\big\{\pi_{\al\,\al+1}(p)\mid p\in U\big\}}$$
For $y\in V^{r.o.(P_{\al+1})}$ we define $\pi''_{\al\,\al+1}$ by
induction on the rank of $y.$ $\pi''_{\al\,\al+1}(\emptyset)=\emptyset$ and 
$$\pi''_{\al\,\al+1}(y)=\big\{(\,\pi''_{\al\,\al+1}(x),\pi'_{\al\,\al+1}(b)\,)\mid (x,b)\in
y\big\}$$

\noindent If $\beta>\alpha+1,$ $\be=\ga+1,$ then for $p\frown\dot{q}\in
P_{\be}$ such that $\dot{q}\neq 1,$
$$\pi_{\al\be}(p\frown\dot{q})=\pi_{\al\ga}(p)\frown \pi''_{\al\ga}(\dot{q})$$
Otherwise, $\pi_{\al\be}(p\frown\dot{q})=\pi_{\al\ga}(p)\frown
1.$
For $p,q\in P_{\al\be},$
$$p\leq_{\al\be}q$$
iff $p\restriction \gamma-\al\leq q\restriction \gamma-\al$ and
$p\restriction \gamma-\al\Vdash p(\gamma-\al) \leq q(\gamma-\al).$ 
If $U\sub P$ is a
regular cut, then
$\pi'_{\al\be}(U)=\overline{\big\{\pi_{\al\be}(p)\mid p\in U\big\}}.$
For $y\in V^{r.o.(P_{\be})}$ we define $\pi''_{\al\be}$ by
induction on the rank of $y.$ $\pi''_{\al\be}(\emptyset)=\emptyset$ and 
$$\pi''_{\al\be}(y)=\big\{(\,\pi''_{\al\be}(x),\pi'_{\al\be}(b)\,)\mid (x,b)\in
y\big\}$$

\noindent If $\be$ is a limit ordinal then if $p\in P_{\be},$ $\pi_{\al\be}=$ the
$\be-\al$ sequence such that for $\zeta<\beta-\al,$ 
$$\pi_{\al\be}(p)\restriction\zeta=\pi_{\al\zeta}(\,p\restriction(\al+\zeta)\,)$$
$\leq_{\al\be}$ is defined in the natural way. If $U\sub P$ is a
regular cut, then
$\pi'_{\al\be}(U)=\overline{\big\{\pi_{\al\be}(p)\mid p\in U\big\}}.$
For $y\in V^{r.o.(P_{\be})}$ we define $\pi''_{\al\be}$ by
induction on the rank of $y.$ $\pi''_{\al\be}(\emptyset)=\emptyset$ and 
$$\pi''_{\al\be}(y)=\big\{(\,\pi''_{\al\be}(x),\pi'_{\al\be}(b)\,)\mid (x,b)\in
y\big\}$$

\noindent Let $\Big({\bf P}_{\al\,Ord}\,,\ {\bf \leq\,}_{\al\,Ord}\Big)$ be the proper class
partial order such that
$${\bf P}_{\al\,Ord} = \bigcup_{\al<\beta\in\,Ord}P_{\al\be}$$
and
$${\bf \leq\,}_{\al\,Ord} = \bigcup_{\al<\beta\in\,Ord}\leq_{\al\be}$$

\vspace{.1in}

\noindent{\bf Theorem 2.} \ Let ${\bf P}$ be a $\psi(x,y)$ definable
$Ord$ iteration. Let
$\alpha\in Ord,$ and let $G$ be a generic subset of $P_{\al}.$ Let
$\varphi(x_1,\ldots,x_n)$ be a formula and $a_1,\ldots, a_n\in
V^{P_{\alpha}}.$ For every $\alpha<\beta\in Ord,$ 

\begin{enumerate}
\item $\pi'_{\alpha\beta}: r.o.(P_{\be})^V\ \longrightarrow\
r.o.(P_{\al\be})^{V[G]}$ is a $\Sigma$ complete Boolean homomorphism 
\item $\pi''_{\alpha\beta}: V^{r.o.(P_{\be})}\ \longrightarrow\ V[G]^{r.o.(P_{\al\be})}$ is onto
\item $\pi'_{\alpha\beta}(\,||\varphi(a_1,\ldots,a_n)||\,)
=||\varphi(\pi''_{\alpha\beta}(a_1),\ldots,\pi''_{\alpha\beta}(a_n))||$
\end{enumerate}

\vspace{.1in}

\noindent We now proceed to prove theorem 2 in the case $\be=\al+1$ with 
the following series of lemmas (3-9) and then prove by induction on $\be$
a further series of 
lemmas (10-14) needed to finish the proof of theorem 2 for a general $\be.$

\vspace{.1in}

\noindent {\bf Lemma 3.} \ Let $p\in P_{\al\,\al+1}.$  There exists $(p_1,\dot{p})\in
P_{\al+1}$ such that 
$$\pi'_{\al\,\al+1}(U_{(p_1,\dot{p})})=U_p$$
\proof Let us denote $\pi'_{\al\,\al+1}$ by $\pi'.$
First note that if $(q_1,\dot{q})\leq (p_1,\dot{p})$ then
$\pi(q_1,\dot{q})\leq \pi(p_1,\dot{p}).$ Let $p\in G$ and $\dot{p}$ a
$P_{\al}$ name such that $i_{G}(\dot{p})=p.$ $\pi(p_1,\dot{p})=p$ and
by the above
$$\pi U_{(p_1,\dot{p})}\sub U_p$$
If $q\leq p$ let $q_1\leq p_1$ and $\dot{q}$ a $P_{\al}$ name such
that $q_1\in G$ and $q_1\Vdash \dot{q}\leq \dot{p}.$ Then
$\pi(q_1,\dot{q})=q$ so $\pi U_{(p_1,\dot{p})}=U_p.$

\vspace{.1in}

\noindent {\bf Lemma 4.} \ Let $(p_1,\dot{p})\in P_{\al+1}$ such that $p_1\in
G.$ For some $p\in P_{\al\,\al+1},$
$$\pi'_{\al\,\al+1}(U_{(p_1,\dot{p})})=U_p$$
\proof Let us denote $\pi_{\al\al+1}$ by $\pi$ and $\pi'_{\al\al+1}$ by $\pi'.$
Let $p=i_G(\dot{p}).$ Since $q_1\in G$ and $(q_1,\dot{q})\leq
(p_1,\dot{p})\ \rightarrow\ \pi(q_1,\dot{q})\leq\pi(p_1,\dot{p})$ we
know 
$$\pi'(U_{(p_1,\dot{p})})\sub U_p$$
Let $q\in U_p$ i.e., $q\leq p.$ Let $q_1\leq p_1$ such that $q_1\Vdash
\dot{q}\leq \dot{p}$ and $q_1\in G.$ Then $\pi(q_1,\dot{q})=q$ and $(q_1,\dot{q})\leq
(p_1,\dot{p}).$ So $\pi'(U_{(p_1,\dot{p})}) = U_p.$

\vspace{.1in}

\noindent {\bf Lemma 5.} \ Let $\pi_{\al\,\al+1}(p_1,\dot{p})=p$ and
$\pi_{\al\,\al+1}(q_1,\dot{q})=q.$ Then $U_{(p_1,\dot{p})}\cap
U_{(q_1,\dot{q})}=\emptyset\ \rightarrow\ U_p\cap U_q=\emptyset.$\\ 
\proof Left to the reader.

\vspace{.1in}

\noindent {\bf Lemma 6.} \ Let $U\sub P$ be a regular cut. Then
$\pi'_{\al\,\al+1}(-U)=-(\pi'_{\al\,\al+1}U).$ \\
\proof Let us denote $\pi_{\al\,\al+1}$ by $\pi$ and
$\pi'_{\al\,\al+1}$ by $\pi'.$ 
Let $(q_1,\dot{q})\in -U.$ So for all $(p_1,\dot{p})\in U,$
$$U_{\pi(p,\dot{p})}\cap U_{\pi(q,\dot{q})}=\emptyset$$
Therefore 
$$\overline{\big\{\pi(p_1,\dot{p})\mid (p_1,\dot{p})\in U\big\}}\cap
U_{\pi(q_1,\dot{q})}=\emptyset$$ 
i.e., $\pi(q,\dot{q})\in -\pi'(U).$ 
Similarly, if 
$$\overline{\big\{\pi(p_1,\dot{p})\mid (p_1,\dot{p})\in U\big\}}\cap
U_{q}=\emptyset$$
then if $r_1\in G$ such that $r_1\Vdash \pi'U\cap U_{q}=\emptyset,$
then for every $(p_1,\dot{p})\in U,$ 
$$U_{(p_1,\dot{p})}\cap U_{(r_1,\dot{q})}=\emptyset$$
so $q=\pi(r,\dot{q})$ and $(r,\dot{q})\in -U.$ (Why? Suppose there
is a $(s_1,\dot{s})\leq (r_1,\dot{q})$ and $(s_1,\dot{s})\leq
(p_1,\dot{p})$ for some $(p_1,\dot{p})\in U.$ 
$$s_1\Vdash \pi'U\cap U_{q}=\emptyset\ \wedge\ s_1\Vdash \dot{s}\in
U_p\cap U_q$$
a contradiction.)

\vspace{.1in}

\noindent {\bf Lemma 7.} \ Let $I$ be an index set and
let $\big\{U_i\mid i\in I\big\}\in V$ be a collection of regular open cuts
in $P_{\al+1}.$ Then
$$\pi'_{\al\,\al+1}(\prod\limits_{i\in I}U_i)=\prod\limits_{i\in
I}\pi'_{\al\,\al+1}U_i$$ 
\proof Let us denote $\pi_{\al\al+1}$ by $\pi$ and $\pi'_{\al\,\al+1}$
by $\pi'.$ Note that
$\prod\limits_{i\in I}U_i=\bigcap\limits_{i\in I}U_i$ for $U_i$
regular cuts of a separative partially ordered set. So
$$\pi'(\prod\limits_{i\in I}U_i)\sub \prod\limits_{i\in I}\pi' U_i$$
is clear. Let $p\in \bigcap\limits_{i\in I}\pi' U_i$ and let
$p_i\frown\dot{p}\in U_i$ such that $\pi(p_i\frown\dot{p})=p.$ (Such a
$p_i\frown\dot{p}$ exists. Let $p_i\in G$ such that $p_i\Vdash
\dot{p}\in \pi'U_i.$ Then $p_i\frown \dot{p}\in U_i.$ For suppose not.
Then there exists $q_i\frown \dot{q}\leq p_i\frown \dot{p}$ such that
$U_{q_i\frown\dot{q}}\cap U_i=\emptyset.$ Let $H$ be a generic subset
of $P_{\al}$ such that $q_i\in H.$ There exists an $r\leq q$ such that
$r=\pi(r_i\frown\dot{r})$ (Here we mean $\pi$ as defined in $V[H]$)
and $r_i\in H$ and $r_i\frown\dot{r}\in U_i.$ Withoutloss of
generality $r_i\leq q_i$ and $r_i\Vdash \dot{r}\leq \dot{q}.$ 
This is a contradiction of $U_{q_i\frown\dot{q}}\cap U_i=\emptyset.$) 
For each $i\in I$ let $U_i'= \big\{r\in P_{\al}\mid
r\frown\dot{p}\in U_i\big\}.$ Each $U_i'$ is a regular cut on
$P_{\al}.$ Why? Suppose $r\in P_{\al}$ and the set of things in $U_i'$ are
dense below $r.$ If $r\not\in U_i'$ then $r\frown\dot{p}\not\in U_i\
\rightarrow$ since $U_i$ is a regular cut there exists $r'\leq r$ and
$\dot{p}'$ such that $r'\frown \dot{p}'\leq r\frown \dot{p}$ and
$$U_{r'\frown\dot{p}'}\cap U_i=\emptyset$$
But $U_i'$ is dense below $r$ so there exists a $r''\leq r'$
such that $r''\frown\dot{p}\in U_i$ which implies
$r''\frown\dot{p}'\in U_i$ a contradiction. Let $G'$ be the
ultrafilter on $r.o.(P_{\al})$ associated with $G.$ For every $i\in I,
\ U_i'\in G'$ which implies since $G'$ is generic that there is an $s\in
\bigcap\limits_{i\in I}U_i'$ such that $s\in G$ which implies
$s\frown\dot{p}\in U_i$ for every $i\in I.$ Since 
$\pi(s\frown\dot{p})=p,$ we have proved
$\pi'_{\al\,\al+1}(\prod\limits_{i\in I}U_i)=\prod\limits_{i\in 
I}\pi'_{\al\,\al+1}U_i.$ 

\vspace{.1in}

\noindent {\bf Lemma 8.} $\pi''_{\al\al+1}:V^{P_{\al+1}}\longrightarrow
V^{P_{\al\al+1}}$ is onto.\\
\proof Let us denote $\pi''_{\al\al+1}$ by
$\pi''.$ By induction on the rank
of $y\in V^{P_{\al\,\al+1}}$ we show that $y$ has an inverse.
$\pi(\emptyset)=\emptyset.$ 
Let $\gamma$ be large enough so that
each element in the domain of $y$ has an inverse in
$V_{\gamma}^{P_{\al+1}}.$ Let $N$ be a name for $y$ in
$V^{P_{\al}}.$  Let $\tilde{y}$ be the name in $V^{P_{\al+1}}$ such
that $z\in \tilde{y}$ iff for some $b\in V^{P_{\al}}$ such that $b$ is a
name for an element of $r.o.(P_{\al\,\al+1}),$ and for some $x\in
V_{\ga}^{P_{\al+1}},$
$$z=\Big(x,\ \Big\{(p_1,\dot{p})\mid p_1\Vdash (\pi''x,b)\in N\
\wedge\ p_1\Vdash \dot{p}\in b\,\Big\}\Big)$$
Then $\pi''(\tilde{y})=y.$ 
(For each $x$ and $b,$ 
$$\Big\{(p_1,\dot{p})\mid p_1\Vdash (\pi''x,b)\in N\
\wedge\ p_1\Vdash \dot{p}\in b\,\Big\}$$ 
is a regular cut. Why? Certainly
$$\Big\{(p_1,\dot{p})\mid p_1\Vdash (\pi''x,b)\in N\
\wedge\ p_1\Vdash \dot{p}\in b\,\Big\}$$ 
is downward closed. Suppose that
$$\Big\{(r_1,\dot{r})\mid r_1\Vdash (\pi''x,b)\in N\
\wedge\ r_1\Vdash \dot{r}\in b\,\Big\}$$ 
is dense below $(p_1,\dot{p}).$
Then $p_1\Vdash (\pi''x,b)\in N$ and $p_1\Vdash b $ is dense below
$\dot{p},$  so $p_1\Vdash \dot{p}\in b$ since $i_{G}(b)$ is a regular cut.)

\vspace{.1in}

\noindent {\bf Lemma 9.} \ If $x_1,\ldots,x_n\in V^{P_{\al+1}}$ and
$\varphi(v_1,\ldots,v_n)$ is a formula, then
$$\pi'_{\al\al+1}(\,||\,\varphi(x_1,\ldots,x_n)\,||\,)=
||\,\varphi(\,\pi''_{\al\,\al+1}x_1,\ldots,\pi''_{\al\,\al+1}x_n\,)\,||$$
\proof First for atomic formulas by induction on
$\Gamma(\rho(x),\rho(y))$ and then by induction on the complexity of
$\varphi(v_1,\ldots,v_n).$ For simplicity we denote $\pi'_{\al\,\al+1}$
as $\pi'$ and $\pi''_{\al\al+1}$
as $\pi''.$ 
$$\pi'(\,||x\in y||\,)\ =\ \pi'\Big(\sum_{t\in dom\,y}||t=x||\bullet
y(t)\ \Big)=\sum_{t\in dom\,y}\pi'(\,||t=x||\,)\bullet\pi'(y(t))\ =$$
$$\sum_{\pi'' t\in dom\,\pi'' y}||\pi'' t=\pi'' x||\bullet \pi''
y(\pi'' t)\ \ =\ \ ||\pi''
x\in \pi'' y||$$

\vspace{.1in}

$$\pi'(\,||x=y||\,)\ =\ \pi'\Big(\prod_{t\in dom\,x}-x(t)+||t\in
y||\ \ \bullet\ \prod_{t\in dom\,y}-y(t)+||t\in x||\ \Big)\ =$$
$$\prod_{t\in dom\,x}\pi'(-x(t))+\pi'(\,||t\in y||\,)\ \ \bullet\ 
\prod_{t\in dom\,y}\pi'(-y(t))+\pi'(\,||t\in x||\,)\ \ =$$
$$\prod_{\pi'' t\in dom\,\pi''\,x}\!\!\!-\pi'' x(\pi'' t)\ +\ ||\pi'' t\in \pi''
y||\ \ \ \bullet\ \prod_{\pi'' t\in dom\,\pi'' y}\!\!\!-\pi'' y(\pi'' t)\ +\ ||\pi''
t\in \pi'' y||\ \ =$$
$$||\pi'' x=\pi'' y||$$

\vspace{.1in}

$$\pi'(\,||\exists x_0\varphi(x_0,x_1,\ldots,x_n)||\,) = \pi'\Big(\sum_{x_0\in
V^{P_{\al+1}}}||\varphi(x_0,x_1,\ldots,x_n)||\ \Big)\ =$$
$$\sum_{x_0\in V^{P_{\al+1}}}\varphi(\pi'' x_0,\pi'' x_1,\ldots,\pi''
x_n)\ =\sum_{y\in V^{P_{\al\,\al+1}}}||\varphi(y,\pi''
x_1,\ldots,\pi'' x_n)||\ =$$
$$||\,\exists x_0\varphi(x_0,\pi'' x_1,\ldots, \pi'' x_n)||\,$$

\vspace{.1in}

\noindent Note that lemma $8$ is used in the next to the last step.
Now we prove a series of lemmas by simultaneous induction
on $\be$ needed to finish the proof of theorem $2.$

\vspace{.1in}

\noindent {\bf Lemma 10.} \ If $p,q\in P_{\be},$ and $p\leq_{\be} q,$ then
$\pi_{\al}p\leq_{\al\be} \pi_{\al}q.$\\
\proof If $\be=\al+1,$ it follows from the definition of $\pi_{\al}.$
If $\be>\al+1,$ it follows from lemma $12$ for ordinals less than $\be.$

\vspace{.1in}

\noindent {\bf Lemma 11.} \ Let $r\in G$ such that $r\frown p_1$ and $r\frown p_2\in
P_{\al+\be}$ such that $r\Vdash \pi_{\al\be}(p_1)=\pi_{\al\be}(p_2).$
Then there is a $s\in G$ such that $s\leq r$ and $s\frown p_1=s\frown
p_2.$ \\
\proof If $\be =\al+1,$ then $r\Vdash p_1\leq
p_2$ and $r\Vdash p_2\leq p_1,$ so $r\frown p_1=r\frown p_2.$ If $\be$
is a limit, then for each $\gamma<\be$ there is a $s_{\gamma}\in G$
such that $s_{\gamma}\frown p_1\restriction\gamma=s_{\gamma}\frown
p_2\restriction\gamma.$ Let $G'$ be the generic ultrafilter on
$r.o.(P_{\al})$ associated with $G.$ If we let $U_{\ga}=\big\{s\in P_{\al}\mid
s\frown p_1\restriction\ga=s\frown p_2\restriction\ga\big\}$ then by
lemma $13$ each
$U_{\ga}$ is a regular cut such that $U_{\ga}\in G'.$ Let $s\in
\bigcap_{\gamma\in\beta}U_{\gamma}\cap G.$
So let $\beta=\gamma+1,\ \gamma>\al.$ $r\Vdash
\pi_{\al\be}p_1=\pi_{\al\be}p_2\ \rightarrow\ r\Vdash
\pi_{\al\ga}p_1=\pi_{\al\ga}p_2\ \rightarrow \ $there is a $t\leq r$
such that $t\in G$ and $t\frown p_1\restriction\ga=t\frown p_2\restriction\ga.$ We
know
$$\pi_{\al\ga}(t\frown p_1\restriction\ga)\Vdash
\pi''_{\al\ga}(p_1(\ga))\leq \pi''_{\al\ga}(p_2(\ga))$$
and vice versa, so by lemma $12$ for ordinals less than $\be,$ there
is a $s'\leq t$ such that $s'\in G$ and 
$$s'\frown p_1\restriction\ga\Vdash p_2(\ga)\leq p_1(\ga)\ \wedge\
p_1(\ga)\leq p_2(\ga)$$
$s'\Vdash
\pi_{\al\ga}(p_1\restriction\ga)=\pi_{\al\ga}(p_2\restriction\ga)$ so
by the induction hypothesis there is an $s\leq s'$ such that $s\in G$
and $s\frown p_1\restriction\ga=s\frown p_2\restriction\ga.$
Therefore, $s\frown p_1=s\frown p_2.$ 

\vspace{.1in} 

\noindent {\bf Lemma 12.} \ Let ${\bf P}$ be a $\psi(x,y)$ definable
$Ord$ iteration. Let $\al<\be\in Ord$ and let $G$ be a generic subset
of $P_{\al}.$ Let $p\frown p_1\in P_{\alpha+\beta}$ such that
$p\in G.$ Let $x_1,\ldots,x_n\in V^{P_{\alpha+\beta}}.$ Then 
$$p\frown p_1\Vdash \varphi(x_1,\ldots,x_n)\ \rightarrow\
\pi_{\al}(p\frown p_1)\Vdash \varphi(\pi''_{\al}
x_1,\ldots,\pi''_{\al} x_n)$$
and if $\pi_{\al\be}(p\frown p_1)\Vdash \varphi(\pi''_{\al\be}
x_1,\ldots,\pi''_{\al\be} x_n)$
then there exists $s\leq p$ such that $s\in G$ and
$$s\frown p_1\Vdash \varphi(x_1,\ldots x_n)$$
\proof By lemma $10,$ $p\leq q\ \rightarrow\ \pi_{\al\be}p\leq \pi_{\al\be}q$
so $U_p\sub ||\varphi(x_1,\ldots,x_n)||\ \rightarrow\
\pi'_{\al\be}(U_p)=U_{\pi_{\al\be}(p)}\sub
||\varphi(\pi''_{\al\be}x_1,\ldots,\pi''_{\al\be}x_n)||.$ Now suppose
$\pi_{\al\be}(p\frown p_1)\in
\pi'_{\al\be}||\varphi(x_1,\ldots,x_n)||.$
Let $q\in G$ such that $q\leq p$ and
$$q\Vdash \pi_{\al\be}(p\frown p_1)\in\pi'||\varphi(x_1,\ldots,x_n)||$$
As in Lemma $7,$ $q\frown p_1\in ||\varphi(x_1,\ldots,x_n)||.$ 

\vspace{.1in}

\noindent {\bf Lemma 13.} \ Let $\be>\al$ and let $p\frown p_1,\
p\frown p_2\in P_{\be}.$ Then
$$U_{p_1,p_2}=\big\{s\in P_{\al}\mid s\frown p_1=s\frown p_2\big\}$$
is a regular cut.\\
\proof By induction on $\be.$ The limit case is easy. So let
$\be=\al+1.$ Let $s\in P_{\al}$ such that for every $s'\leq s$ there is
a $t\leq s'$ such that $t\frown p_1=t\frown p_2.$ Then $t\Vdash
p_1=p_2.$ Since the set of $t$ below $s$ which force $p_1=p_2$ is
dense below $s,$ $s\Vdash p_1=p_2.$ So $s\frown p_1=s\frown p_2$ i.e.,
$s\in U_{p_1,p_2}.$ So let $\be=\ga+1$ with $\ga>\al.$ Let $\big\{t\in
P_{\al}\mid t\frown p_1=t\frown p_2\big\}$ be dense below $s.$ By the
induction hypothesis, $s\frown p_1\restriction\ga=s\frown
p_2\restriction\ga.$ If $s'\frown p'\leq s\frown p_1,$ let $t\leq s'$
such that $t\frown p_1\restriction\ga\Vdash p_1(\ga)=p_2(\ga).$ Then
$t\frown p'\restriction\ga\Vdash p_1(\ga)=p_2(\ga)$ i.e., the set of
things below $s\frown p_1\restriction\ga$ forcing $p_1(\ga)=p_2(\ga)$
is dense below $s\frown p_1(\ga),$ so $s\frown p_1\restriction\ga\Vdash
p_1(\ga)=p_2(\ga).$ So $s\frown p_1=s\frown p_2.$ 

\vspace{.1in}

\noindent {\bf Lemma 14.} \ Let $p\frown p_1,\ p\frown p_2\in P_{\be}$ such
that $p\in P_{\al}.$ Let $\dot{G}$ be the canonical name for a generic
subset of $P_{\al}$ and let $r\leq p.$ Suppose
$$r\Vdash p\in \dot{G}\ \wedge\ \pi_{\al}(p\frown p_1)\leq
\pi_{\al}(p\frown q_1)$$
Then $r\frown p_1\leq r\frown q_1.$ \\
\proof By induction on $\be$ where $p\frown p_1, p\frown q_1\in
P_{\be}.$ If $\be=\al+1$ then it follows from the definition of $\leq$
on $P_{\al+1}.$ The limit case is also easy. So let $\be=\ga+1,\
\ga>\al.$ By induction we can assume $r\frown p_1\restriction\ga\leq
r\frown q_1\restriction\ga.$ We must show
$$r\frown p_1\restriction\ga\Vdash p_1(\ga)\leq q_2(\ga)$$
We have by assumption
$$r\Vdash \Big(\ \pi_{\al\ga}(p_1\restriction\ga)\ \Vdash\
\pi''_{\al\ga}(p_1(\ga))\ \leq
\ \pi''_{\al\ga}(q_1(\ga))\ \Big)$$
Let $r'\leq r.$ Let $H$ be a generic
subset of $P_{\al}$ such that $r'\in H.$ By lemma $12$ there is an
$s\in H$ such that $s\leq r'$ and
$$s\frown p_1\restriction\ga\Vdash p_1(\ga)\leq q_1(\ga)$$
We have just shown that the set of things in $P_{\ga}$ forcing
$p_1(\ga)\leq q_1(\ga)$ is dense below $r\frown p_1\restriction \ga.$ So
$r\frown p_1\restriction\ga\Vdash p_1(\ga)\leq q_1(\ga).$ 

\vspace{.1in} 

\noindent \proof (theorem $2$) The proof of theorem $2$ now follows with
slight modifications over the case $\be=\al+1,$ using 
lemmas 11 and 13. The details are left to the reader. 

\vspace{.1in}

\noindent {\bf Corollary 15.} \ Let ${\bf P}$ be a $\psi(x,y)$ definable
$Ord$ iteration. Let
$\alpha\in Ord,$ and let $G$ be a generic subset of $P_{\al}.$
Then ${\bf P}_{\al\,Ord}$ is a $\psi(x,y+\al)$
definable iteration in $V[G].$\\
\proof  Note that a $\psi(x,y+\al)$ definable $Ord$ iteration is
unique up to the choice of parameter in $\psi(x,y),$ so without loss
of generality we can speak of the $\psi(x,y+\al)$ definable $Ord$
iteration ${\bf P}.$ By induction on $\be>\al,$ we prove that every
$\beta-\al$ sequence in ${\bf P}_{\al\,Ord}$ is in ${\bf P}$ and every
$\beta-\al$ sequence in ${\bf P}$ is in ${\bf P}_{\al\,Ord}.$ The
details are left to the reader, but one uses theorem $2$ together with
lemma $12.$ 

\vspace{.1in}

\noindent {\bf Theorem 16.} \ Let ${\bf P}$ be a $\psi(x,y)$ definable
$Ord$ iteration. Let 
$\al<Ord$ and let $\pi=\pi_{\al}.$ Let ${\bf G}$ be a $V$ generic
subclass of ${\bf P}.$ Let $G={\bf G}\cap P_{\al}$ and ${\bf
H}=\pi[{\bf G}].$ Then
\begin{enumerate}
\item $G$ is a $V$ generic subset of $P_{\al}$
\item ${\bf H}$ is a $V[G]$ generic subclass of ${\bf P}_{\al\,Ord}$ 
\item $V[{\bf G}]=V[G][{\bf H}]$
\end{enumerate}
\proof The proof of $1$ and the proof that ${\bf H}$ is a filter is
left to the reader. To finish $2,$ let 
$\theta(x,a)$ be a formula with $a\in V[G]$ such that $\theta(x,a)$
defines a dense subclass of ${\bf P}_{\al\,Ord}.$ Let $\dot{a}$ be a
name for $a.$ Let $p\in G$ such that $p\Vdash \theta(x,\dot{a})$
defines a dense subclass $D'$ of ${\bf P}_{\al\,Ord}.$ Let $p\frown p_1\in
{\bf G}.$ Let $q\frown q_1\leq p\frown p_1.$ Let $r\leq q$ such that
for some $r\frown r_1\in {\bf P},$
$$r\Vdash \pi r_1\leq \pi q\ \wedge\ \varphi(\pi r_1,\dot{a})$$
By lemma $14,$ $r\frown r_1\leq q\frown q_1.$ So the set of $r\frown r_1$
such that $r\Vdash \varphi(\pi r_1,\dot{a})$ is dense below $p\frown
p_1.$ So for some $r\frown r_1\in {\bf G},\  r\Vdash \varphi(\pi
r_1,\dot{a}),$ i.e., $\pi r_1\in \pi[{\bf G}]\cap D'={\bf H}\cap D'.$
So ${\bf H}$ is a $V[G]$ generic subclass of ${\bf P}_{\al\,Ord}.$ 
To prove $3,$  we show by induction on the rank of $x\in V^{\bf P}$ that
$i_{\bf G}(x)=i_{\bf H}(\pi x).$ By definition, $\pi \emptyset=\emptyset,$ so 
$i_{\bf G}(\emptyset)=i_{\bf H}(\emptyset)=\emptyset.$ Let $y$ have
rank $\ga+1$ and suppose the theorem is true for all names with rank
$\leq \ga.$ Then $i_{\bf G}(y)=\big\{i_{\bf G}(x)\mid y(x)\cap {\bf
G}\neq\emptyset\big\}=\big\{i_{\bf H}(\pi x)\mid \pi y(x)\cap \pi[{\bf
G}]\neq\emptyset\big\}=i_{\bf H}(\pi y).$ (Why? Let $S\subseteq {\bf P}$ be
a regular cut of ${\bf P}.$ ${\bf G}\cap S\neq\emptyset\ \rightarrow\
{\bf H}\cap \pi[S]\neq\emptyset.$ On the other hand if ${\bf H}\cap
\pi[S]\neq \emptyset,$ let $p\frown p_1\in {\bf G}$ and $s\frown
 s_1\in S$ such that $\pi(p\frown p_1)=\pi(s\frown s_1).$ Let $r\in G$
such that $r\leq s\ \wedge\ r\leq p\ \wedge\ r\Vdash \pi(p\frown
p_1)=\pi(s\frown s_1)$ which
implies by lemma $11$ that there exists a $s'\in G$ such that $s'\frown
s_1=s'\frown p_1.$ Since $S$ is a cut $s'\frown s_1\in S.$ Since
$s'\in G$ there exists $p_1'$ such that $s'\frown p_1'\in {\bf G}.$
$s'\frown p_1'\in {\bf G}$ and $p\frown p_1\in {\bf G}$ implies there is
a $t\frown t'\in {\bf G}$ such that $t\frown t_1'\leq s'\frown p_1'$
and $t\frown t_1'\leq p\frown p_1.$ $t\frown t_1'\leq s'\frown p_1\
\rightarrow\  s'\frown p_1\in {\bf G}.$ )

\vspace{.1in}

\noindent {\bf Lemma 17.} \ Let ${\bf P}$ be a $\psi(x,y)$ definable
$Ord$ iteration. Suppose $P_1$ is $\aleph_{\al}$ closed and that for every
$\be\in Ord$ and $p\in P_{\be},$ if $p\Vdash \exists !P(\psi(P,\al)$
then $p\Vdash P$ is $\aleph_{\al}$ closed. Then ${\bf P}$
is $\aleph_{\al}$ closed. \\
\proof Left to the reader.

\vspace{.1in}

\noindent {\bf Lemma 18.} \ Let ${\bf P}$ be a class forcing such that ${\bf
P}$ is $\aleph_{\al}$ closed. Let ${\bf G}\sub {\bf P}$ be a generic
subclass of ${\bf P}.$ Then if $f$ is a function from $\aleph_{\al}$
into $V$ and $f\in V[{\bf G}],$ then $f\in V.$ \\
\proof Standard.

\vspace{.1in}

\begin{center}
{\bf A CLASS FORCING DEMONSTRATING THE CONSISTENCY OF $ZFC\
+\ CIFS$}
\end{center}

\vspace{.1in}

In this section given an arbitrary finite list
$\{\psi_1,\ldots,\psi_n\}$  of formulas we define in $L$ a class
partial ordering ${\bf P}$ definable in $L$ such that ${\bf P}$ is a
$\psi(x,y)$ definable iteration for some formula $\psi(x,y)$ (an
unwieldy combination of the $\{\psi_1,\ldots,\psi_n\}$) and
$$L^{\bf P}\models \ \ ZFC\ +\ CIFS\restriction\{\psi_1,\ldots, \psi_n\}$$

\vspace{.1in}

\noindent {\bf Definition 3.} \ Let $\aleph_{\al}$ be a regular cardinal
and $X$ a set. 
$$Col(X,\aleph_{\al})$$ 
is the partial order consisting
of all injections of cardinality less
than $\aleph_{\al}$ from $X$ into $\aleph_{\al}.$

\vspace{.1in}

\noindent {\bf Definition 4.} \ We define the class partial ordering
${\bf P}=\bigcup\limits_{\al\in 
Ord-\{0\}}P_{\al}$ and ${\bf\leq}=\bigcup\limits_{\al\in
Ord-\{0\}}\leq_{\al}$ where $P_{\al}$ consists of $\al$ sequences and
$P_{\al}$ and $\leq_{\al}$ are defined
by induction on $\al$ as follows:
$$P_1=\Big\{(0,(p_1,\ldots,p_n)\,)\mid 
p_i=1\ \vee\ \exists !X(\,\psi_i(X)\
\wedge\ p_i\in Col(X,\aleph_{\al})\,)\Big\}$$

\vspace{.15in}

\noindent $\leq_1\ =\ \Big\{(\,(0,(p_1,\ldots,p_n)\,),\ (0,(q_1,\ldots,q_n)\,)\,)\mid q_i=1\ \vee$

\begin{flushright}
$\exists !X(\,\psi_i(X)\ \wedge\ p_i\in Col(X,\aleph_{\al})\ \wedge\
q_i\in Col(X,\aleph_{\al})\ \wedge\ 
p_i\leq q_i\,)\Big\}$
\end{flushright}

\vspace{.1in}

\noindent If $\al$ is a limit ordinal then $P_{\al}=\big\{p\mid p$ is an $\al$
sequence and $\forall\beta\neq 0,\ p\restriction\beta\in P_{\be}\big\}$
and $\leq_{\al}=\big\{(p,q)\mid p,q\in P_{\al}\ \wedge
\ p\restriction\beta\leq_{\be}q\restriction\beta\ \ \forall\be\in
\al-\{0\}\big\}.$ 
If $\al$ is a limit ordinal, then

\vspace{.15in}
 
\noindent $P_{\al+1}=
\Big\{p\frown (\dot{p}_1,\ldots,\dot{p}_n)\mid p\in P_{\al}\ \wedge\
p_i =1\ \ \vee\ \  \dot{p}_i\hbox{ is a }P_{\al}\hbox{ name }\ \wedge$
\begin{flushright}
$p\Vdash \aleph_{\al}\hbox{ is regular and }X\hbox{ is definable by }\psi_i\hbox{ in
}L(V_{\omega+\al})\ \wedge\ p_i\in Col(X,\aleph_{\al})\Big\}$ 
\end{flushright}

$$\leq_{\al+1}=\Big\{(p\frown (\dot{p}_1,\ldots,\dot{p}_n),q\frown
(\dot{q}_1,\ldots,\dot{q}_n))\mid 
p\leq_{\al}q\ \wedge\ q_i=1\ \vee\ p\Vdash p_i\leq q_i\Big\}$$

\noindent If $\al$ is a successor, then 

\vspace{.15in}

\noindent $P_{\al+1}=\Big\{p\frown (\dot{p}_0,\dot{p}_1,\ldots,\dot{p}_n)\mid p\in
P_{\al}\ \wedge$
\begin{center}
$p\Vdash \dot{p}_0 \hbox{ is a bijection from an ordinal less than
}\aleph_{\al}\hbox{ into }V_{\omega+\al}$
\end{center}
\begin{center}
$\wedge\ \ 1\leq i\leq n\ \rightarrow\ \dot{p}_i=1\ \ \vee$
\end{center}
\begin{flushright}
$p\Vdash \exists !X(\,L(V_{\omega+\al})\models \psi_i(X)\
\wedge\ \dot{p}_i\in Col(X,\aleph_{\al})\,)\,\Big\}$
\end{flushright}

\vspace{.1in}

\noindent $\leq_{\al+1}$ is defined as in the case $\al$ is a limit. 

\vspace{.1in}

\noindent {\bf Lemma 19.} \ Let ${\bf P}$ be as above. Let ${\bf G}$ be a
$L$ generic subclass of ${\bf P}$ and if $\al\in Ord,$ let
$G_{\al}={\bf G}\cap P_{\al}.$ Then
\begin{enumerate}
\item $L[{\bf G}]\models ZFC$
\item $V_{\omega+\al}^{L[G_{\al}]}=V_{\omega+\al}^{L[{\bf G}]}$
\item $L[G_{\al+1}]\models |V_{\omega+\al}|=\aleph_{\al}$
\item $\aleph_{\al}^{L[G_{\al}]}=\aleph_{\al}^{L[{\bf G}]}$
\end{enumerate}
\proof The proof that $L[{\bf G}]\models ZFC$ uses corollary
$15,$ theorem $16,$ and lemmas $17$ and $18$ and is similar to the
proof that the class forcing extension used to prove Easton's theorem
satisfies $ZFC.$ The rest of the proof of
lemma $19$ is by induction on $\alpha.$ If $\alpha=\be+1,$ then ${\bf P}_{\al\,Ord}$
is $\aleph_{\beta}$ closed and $L[{\bf G}]=L[G_{\al}][{\bf H}]$ where
${\bf H}$ is a $L[G_{\al}]$ generic subclass of ${\bf P}_{\al\,Ord},$
so by lemma $18,$ $V_{\omega+\al}^{L[G_{\al}]}=V_{\omega+\al}^{L[{\bf
G}]}$ and $\aleph_{\al}^{L[G_{\al}]}=\aleph_{\al}^{L[{\bf G}]}.$
$L[G_{\al+1}]\models |V_{\omega+\al}|=\aleph_{\al}$ by the definition of
${\bf P}$ and the fact that $V_{\omega+\al}^{L[G_{\al}]}=V_{\omega+\al}^{L[{\bf G}]}.$
If $\al$ is a limit it follows from the induction hypothesis.

\vspace{.1in}

\noindent {\bf Lemma 20.} \ Let ${\bf P}$ be as above. Let ${\bf G}$ be a
$L$ generic subclass of ${\bf P}$ and if $\al\in Ord,$ let
$G_{\al}={\bf G}\cap P_{\al}.$ If $L[{\bf G}]\models
\aleph_{\al}$ is regular$\ \wedge\ \big(\,L(V_{\omega+\al})\models \exists !X(\psi_i(X))\,\big)$
 then $L[{\bf G}]\models |X|\leq \aleph_{\al}.$\\
\proof  If $L[{\bf G}]\models
\big(\,L(V_{\omega+\al})\models \exists !X(\psi_i(X))\,\big),$ then
since $V_{\omega+\al}^{L[G_{\al}]}=V_{\omega+\al}^{L[{\bf G}]}$ we
also have, $L[G_{\al}]\models 
\big(\,L(V_{\omega+\al})\models \exists !X(\psi_i(X))\,\big).$ Since
$\aleph_{\al}^{L[G_{\al}]}=\aleph_{\al}^{L[{\bf G}]}$ and $L[{\bf
G}]=L[G_{\al}][{\bf H}]$ where 
${\bf H}$ is a $L[G_{\al}]$ generic subclass of ${\bf P}_{\al Ord},$
${\bf H}\restriction ({\bf P}_{\al Ord})_1$ witnesses the existence of
an bijection from $X$ onto $\aleph_{\al}.$

\pagebreak

\begin{center}
REFERENCES
\end{center}

\begin{enumerate}

\item {[Jech1] T. Jech, {\em Multiple Forcing}, Cambridge
University Press.}

\item {[Jech2] T. Jech, {\em Set Theory}, Academic Press.}

\end{enumerate}

\end{document}